\documentclass[11pt, a4paper]{article}
\usepackage[a4paper,margin=3cm,left=3cm,right=3cm]{geometry}
\usepackage[english]{babel}

\usepackage{amsmath}
\usepackage{amssymb}
\usepackage{bm}
\usepackage{bbm}
\usepackage{graphicx}
\usepackage[colorlinks=true, allcolors=black]{hyperref}
\usepackage{authblk}

\usepackage[caption=false]{subfig}

\usepackage{natbib}
\usepackage{mathtools}

\DeclarePairedDelimiterX{\infdivx}[2]{(}{)}{
  #1\;\delimsize\|\;#2%
}

\newcommand{\ud}{\text{d}}

\usepackage{libertine}
\usepackage[libertine]{newtxmath}

\newtheorem{remark}{Remark}[section]

\usepackage{algorithm}
\usepackage{algpseudocode}

\title{Stochastic parameterisation: the importance of nonlocality and memory}
\author{Martin T. Brolly}
\affil{School of Mathematics and Maxwell Institute for Mathematical Sciences, University of Edinburgh,
King’s Buildings, Edinburgh EH9 3FD, UK}

\begin{document}

\maketitle
\begin{abstract}
    Stochastic parameterisations deployed in models of the Earth system frequently invoke locality assumptions such as Markovianity or spatial locality. This work highlights the impact of such assumptions on predictive performance. Both in terms of short-term forecasting and the representation of long-term statistics, we find locality assumptions to be detrimental in idealised experiments. We show, however, that judicious choice of Markovian parameterisation can mitigate errors due to assuming Markovianity. We propose a simple modification to standard Markovian parameterisations, which yields significant improvements in predictive skill while reducing computational cost. We further note a divergence between configurations of a parameterisation which perform best in short-term prediction and those which best represent time-invariant statistics.
\end{abstract}

\section{Introduction}

In the simulation of physical systems, such as the atmosphere or ocean, it is often impossible to resolve dynamics on the full range of relevant scales at once. Limited by computational cost, numerical models covering planetary scales can offer only a truncated view of such systems. Parameterisations are the attempts made to mitigate the errors introduced by failing to resolve certain dynamics, typically by supplementing the equations solved in our models with additional terms. The classical approach to parameterisation is to introduce terms which diagnose approximately the influence of the unresolved dynamics on the resolved dynamics. Stochastic parameterisation arises when it is inappropriate to assume that a deterministic relationship holds. Instead the probabilistic relationship between the unresolved and resolved dynamics is modelled with stochastic terms.

Stochastic parameterisations promise benefits for both `weather' and `climate' modelling, interpreted broadly as dynamical modelling with the aim of capturing finite-time and long-time average behaviour, respectively. In the weather paradigm the primary advantage of stochastic parameterisations is that they allow to quantify uncertainty due to model error in ensemble simulation. Ensembles are commonly used to quantify the effect of uncertainty in initial conditions, however, without additionally capturing uncertainty due to model error, ensemble forecasts can misrepresent the overall uncertainty. Note that in this work we focus solely on model uncertainty. Aside from quantifying model uncertainty, ensembles can be used to produce improved point-estimates, i.e. the ensemble mean of a model with stochastic parameterisations may provide a better prediction than that of a model with deterministic parameterisations~\citep{gneiting2005weather, leutbecher2008ensemble}. In the climate paradigm one hopes that the added variability introduced by stochastic parameterisations leads to a better representation of equilibrium statistics. Indeed there is evidence for all of the above. Stochastic parameterisations have seen significant successful use over the past three decades~\citep{berner2017}, both in research~\citep{wilks2005, crommelin2008, kwasniok2012, kwasniok2014, arnold2013, mana2014, chorin2015, grooms2016, guillaumin2021, levine2022} and in operational weather forecasting~\citep{toth1997, buizza1999, berner2009, palmer2009}.

In constructing stochastic parameterisations it is common to make simplifying assumptions about the relationship between resolved and unresolved variables. E.g. that, given the resolved state, the model error is uncorrelated or Markovian in time and/or space, or that it is Gaussian. This article discusses the impact of such assumptions on the performance of stochastic parameterisations.

The article is structured as follows. In Section~\ref{sec: math} we describe the problem in a general mathematical notation. We also discuss the simplifying assumptions mentioned above in detail. In Section~\ref{sec: 63} we present idealised numerical experiments with the Lorenz '63 model~\citep{lorenz1963}, wherein the original system is forced with artificial model error generated with non-Markovian and spatially-correlated stochastic processes. We show that Markovian and spatially-local parameterisations can fail to reproduce the behaviour of the forced system. On the other hand we show that judicious choice of Markovian approximation can alleviate this error. In Section~\ref{sec: 96} we explore the same issues in the two-scale Lorenz '96 system~\citep{lorenz1996}, wherein the model error arises from neglecting the small-scale variables. We find that modelling spatial correlation in the model error, given the large-scale variables, is critical. We also propose a method for obtaining improved Markovian parameterisations at no additional cost. In Section~\ref{sec: conc} we conclude and discuss the implications of these results for weather and climate modelling in particular.

\section{Mathematical formulation}\label{sec: math}

Here we introduce the problem of parameterisation using a generic notation. Suppose the dynamical system of interest has state $\bm{x}_n$, which can be decomposed as $\bm{x}_n\eqqcolon(\bm{x}_n^{(1)},\,\bm{x}_n^{(2)})$, where $\bm{x}_n^{(1)}$ is a reduced state to be resolved explicitly in a numerical model, and $\bm{x}_n^{(2)}$ is a component to be neglected. 
In the case of turbulence closure, $\bm{x}_n$ would represent the state at some reference high resolution, $\bm{x}_n^{(1)}$ a coarse-grained representation, and $\bm{x}_n^{(2)}$ small-scale fluctuations about the coarse-grained state.
Alternatively, $\bm{x}_n^{(1)}$ could represent the ocean and atmosphere components of an Earth system model, and $\bm{x}_n^{(2)}$ a sea ice or land surface component. In any case, if the full dynamics is given by
\begin{align}\label{eq:system_full}
    \bm{x}_{n+1}=\bm{\Psi}(\bm{x}_n),
\end{align}
then we may equivalently write
\begin{subequations}\label{eq:system_new}
    \begin{align}
        \bm{x}^{(1)}_{n+1} &= \bm{\Psi}_1(\bm{x}_n^{(1)},\,\bm{x}_n^{(2)}),\\
        \bm{x}^{(2)}_{n+1} &= \bm{\Psi}_2(\bm{x}_n^{(1)},\,\bm{x}_n^{(2)}),
    \end{align}
\end{subequations}
where $\bm{\Psi}_1$ and $\bm{\Psi}_2$ denote, separately, the dynamics of $\bm{x}_n^{(1)}$ and $\bm{x}_n^{(2)}$, respectively. Moreover, suppose we have access to an approximate model for $\bm{x}_n^{(1)}$, $\bm{\Psi}_0(\bm{x}_n^{(1)})\approx\bm{\Psi}_1(\bm{x}_n^{(1)},\,\bm{x}_n^{(2)})$, so that
\begin{align}\label{eq: psi0}
    \bm{x}_{n+1}^{(1)}\approx \bm{\Psi}_0(\bm{x}_n^{(1)}).
\end{align}
In practice, an approximate model can usually be obtained by neglecting terms in $\bm{\Psi}_1$ which couple $\bm{x}_n^{(1)}$ to $\bm{x}_n^{(2)}$. In turbulence closure, this amounts to employing the same numerical method (i.e. partial differential equation solver) at coarse resolution.

Given an approximate model as in~\eqref{eq: psi0}, a model error process
\begin{align}\label{eq: m}
    \bm{m}(\bm{x}_n^{(1)}, \bm{x}_n^{(2)}) \coloneqq \bm{\Psi}_1(\bm{x}_n^{(1)},\,\bm{x}_n^{(2)}) - \bm{\Psi}_0(\bm{x}_n^{(1)})
\end{align}
is naturally defined. In the example of turbulence closure, $\bm{m}$ typically corresponds to Reynolds stresses. By definition of $\bm{m}$, we can rewrite~\eqref{eq:system_new} as
\begin{subequations}\label{eq:corrected_system}
    \begin{align}\label{eq:corrected_system_a}
        \bm{x}^{(1)}_{n+1} &= \bm{\Psi}_0(\bm{x}^{(1)}_n) + \bm{m}(\bm{x}^{(1)}_n,\,\bm{x}^{(2)}_n),\\
        \label{eq:corrected_system_b}
        \bm{x}^{(2)}_{n+1} &= \bm{\Psi}_2(\bm{x}_n^{(1)},\,\bm{x}_n^{(2)}).
    \end{align}
\end{subequations}
In fact, by simple substitution, one can show that for any $n\geq1$, there is a function $\bm{m}_n$ such that
\begin{align}\label{eq:corrected_system_history}
    \bm{x}^{(1)}_{n+1}=\bm{\Psi}_0(\bm{x}^{(1)}_n) + \bm{m}_n\left(\{\bm{x}^{(1)}_s\}_{s=0}^n,\,\bm{x}^{(2)}_0\right).
\end{align}
While~\eqref{eq:corrected_system_history} avoids explicit dependence on $\bm{x}_n^{(2)}$, this comes at the cost of dependence on the full history of the reduced state $\{\bm{x}^{(1)}_s\}_{s=0}^n$.

The goal, in any form of parameterisation, is to model the error process, such that $\bm{m}_n$ may be replaced in~\eqref{eq:corrected_system_history} by more tractable terms. That is, given an estimate or approximation $\widehat{\bm{m}}_n$ of $\bm{m}(\bm{x}_n^{(1)}, \bm{x}_n^{(2)})$, one constructs a \textit{parameterised} model
\begin{align}\label{eq:param}
    \widehat{\bm{x}}_{n+1}^{(1)} = \bm{\Psi}_0(\widehat{\bm{x}}_n) + \widehat{\bm{m}}_n.
\end{align}
Conventional parameterisations employ a deterministic approximation of the form $\widehat{\bm{m}}_n=\bm{f}\left(\bm{x}_n^{(1)}\right)$ for some function $\bm{f}$, and thereby neglect dependence on the other arguments of $\bm{m}_n$. The form of $\bm{f}$ is often motivated by physical intuition or asymptotic theories, and simple forms are favoured for tractability of analysis and ease of implementation. Classical examples are the Gent--McWilliams parameterisation~\citep{gent1990}, used to represent the effects of mesoscale eddies in large-scale ocean models, or the Smagorinsky model of eddy viscosity~\citep{smagorinsky1963}. A more modern example is the parameterisation of sub-grid momentum forcing proposed by~\citet{zanna2020} by way of an equation discovery methodology. However, even for particularly simple choices of $\bm{f}$, it is challenging to quantify how the error introduced by these parameterisations propagates into finite-time forecasts and stationary statistics. It is also clear that neglecting the other arguments of $\bm{m}_n$ introduces uncertainty.

Stochastic parameterisations attempt to represent this uncertainty explicitly by modelling $\bm{m}$ probabilistically. It is illustrative to first think of the joint evolution of $\bm{x}_n^{(1)}$ and $\bm{m}_n$ as a stochastic process denoted $(\bm{X}_n^{(1)},\,\bm{M}_n)_{n\in\mathbb{Z}}$ satisfying
\begin{align}\label{eq:approx_dynamics}
    \bm{X}^{(1)}_{n+1} = \bm{\Psi}_0(\bm{X}_n^{(1)}) + \bm{M}_n.
\end{align}
Given knowledge of the full history of the system, this stochastic process is degenerate in the sense that its future evolution is determined. Formally, we may write
\begin{align}\label{eq:dirac_M}
    \bm{M}_n\mid\{\bm{X}_s^{(1)}\}_{s=0}^n,\, \bm{x}_0^{(2)} \sim \delta_{\bm{m}_n\left(\{\bm{x}^{(1)}_s\}_{s=0}^n,\,\bm{x}^{(2)}_0\right)},
\end{align}
and correspondingly
\begin{align}\label{eq:X_dirac}
    \bm{X}^{(1)}_{n+1}\mid\{\bm{X}_s^{(1)}\}_{s=0}^n,\, \bm{x}_0^{(2)} \sim \delta_{\bm{x}^{(1)}_{n+1}},
\end{align}
where $\delta_{\bm{a}}$ denotes a Dirac (probability) measure, assigning probability $1$ to the value $\bm{a}$. Under suitable assumptions, namely that both~\eqref{eq:system_full} and~\eqref{eq: psi0} are ergodic, the joint process $(\bm{X}_n^{(1)},\,\bm{M}_n)_{n\in\mathbb{Z}}$ is stationary. This means that we can meaningfully marginalise with respect to some of the conditioning variables in~\eqref{eq:dirac_M} and consider, for example, the conditional distribution $\bm{M}_n\mid\bm{X}^{(1)}_n$. Crucially, this conditional distribution is not degenerate: $\bm{X}^{(1)}_n$ alone does not determine $\bm{M}_n$. More generally, neglecting dependence on any argument of $\bm{m}_n$ introduces uncertainty.

A common approach to stochastic parameterisation is to assume the process $(\bm{X}_n^{(1)},\,\bm{M}_n)_{n\in\mathbb{Z}}$ is Markovian. Under this assumption a parameterisation is then obtained in a natural way, by drawing realisations $\widehat{\bm{M}}_n$ from the conditional distribution $\bm{M}_n\mid \bm{X}^{(1)}_n,\, \bm{M}_{n-1}$, or simply $\bm{M}_n\mid \bm{X}^{(1)}_n$ and iterating as in~\eqref{eq:param}. While Markovian parameterisations aim to represent the uncertainty introduced by neglecting dependence on the systems history, the Markovian assumption is generally invalid, and hence introduces some error. I.e. the parameterised process $(\widehat{\bm{X}}^{(1)}_n,\,\widehat{\bm{M}}_n)$ may differ significantly from the true process $(\bm{X}_n^{(1)},\,\bm{M}_n)$. Nevertheless Markovian parameterisations remain attractive, since the construction of non-Markovian parameterisations is difficult and costly. While there has been work towards constructing non-Markovian parameterisations, both using approximations guided by theory~\citep{wouters2012, kondrashov2015, vissio2018, gutierrez2021} and by purely data-driven approaches~\citep{crommelin2021, parthipan2023}, these studies have been so far restricted to applications in low dimensional models and/or rely on assumptions of scale separation that are not generally valid. Only Markovian parameterisations are used in operational models. Examples are the popular ``Stochastically Perturbed Parameterisation Tendencies'' (SPPT)~\citep{buizza1999} and ``Stochastically Perturbed Parameters'' (SPP)~\citep{christensen2015} schemes, which are used by meteorological services internationally~\citep{leutbecher2017, christensen2022}.

Additionally, assumptions are often made regarding the spatial correlation structure of the joint process $(\bm{X}_n^{(1)},\,\bm{M}_n)_{n\in\mathbb{Z}}$. For example, in SPPT, two-dimensional random fields are generated which perturb existing deterministic parameterisations multiplicatively, and the same perturbations are applied at each vertical level. While SPPT can be tuned empirically, on the basis of coarse-grained high-resolution simulations (see~\citet{christensen2020} for a detailed discussion), the accuracy that can be attained is ultimately limited by the (quite restrictive) assumptions underlying its construction.

Through a series of numerical experiments we examine the impact of (i) the Markovianity assumption and (ii) spatial autocorrelation (or nonlocality) on the performance of stochastic parameterisations. We also highlight that the natural choice of Markovian parameterisation is generally suboptimal, and introduce a simple method to generate a family of alternative Markovian approximations, over which a modeller can optimise.

\section{Experiments with the Lorenz '63 system}\label{sec: 63}
The experiments in this section make use of the Lorenz '63 system~\citep{lorenz1963}, usually written
\begin{subequations}\label{eq: l63}
    \begin{align}
        \dot{x} &= \sigma(y - x),\\
        \dot{y} &= x(\rho-z) - y,\\
        \dot{z} &= xy - \beta y,
    \end{align}
\end{subequations}
with the classical parameter values $\{\sigma=10,\,\rho=28,\,\beta=8/3\}$. In particular, we consider the discrete-time system obtained through applying the fourth-order Runge--Kutta scheme augmented with artificial model error processes. To align with our notation we write
\begin{align}
    X^{(1)}_{n+1} = \Psi_0\left(X_n^{(1)}\right) + M_n,
\end{align}
letting $X^{(1)}=(x,\,y,\,z)$, $\Psi_0$ denote the map of the discretised Lorenz '63 system, and $M_n$ a generic error process. In sections~\ref{subsec: non-markovian} and~\ref{subsec: nonlocal} we will consider $M_n$ given by non-Markovian and spatially-correlated stochastic processes, respectively. We will then assess the effectiveness of Markovian and spatially local parameterisations as approximations to the true error process. This setup, though artificial, is attractive because it allows us complete control of the properties of the error process. In particular, we can avoid the conflation of errors due to multiple invalid assumptions by treating each in isolation with a targeted experiment.

\subsection{Non-Markovian additive AR-2 forcing}\label{subsec: non-markovian}

Here we consider non-Markovian model error given by the AR$(2)$ process
\begin{align}\label{eq:ar2}
    M_n = \varphi_1 M_{n-1}+\varphi_2 M_{n-2} + \bm{\varepsilon}_n,\quad \bm{\varepsilon}_n\overset{\mathrm{iid}}{\sim}\mathcal{N}(0,\, \sigma_{\varepsilon}^2),
\end{align}
applied independently to each component $(x,\,y,\,z)$ of~\eqref{eq: l63}. The AR$(2)$ process is chosen to highlight that non-Markovianity need not be severe or complex to have important effects. Indeed the AR$(2)$ process is second-order Markovian and we deliberately choose parameter values which should be considered fairly mild, both with respect to non-Markovianity and variance: $\varphi_1=0.45$, $\varphi_2=0.5$, and $\sigma_{\varepsilon}^2=1.425\times10^{-5}$, such that $\mathrm{Var}(M_n)=10^{-4}$.

There is a natural Markovian approximation to the AR$(2)$ process given by the AR$(1)$ process
\begin{align}\label{eq:ar1}
    M_n = \varphi_M M_{n-1} + \epsilon _n,
\end{align}
where $\varphi_M = \frac{\varphi_1}{1-\varphi_2}$ and $\epsilon _n\overset{\text{iid}}{\sim}\mathcal{N}(0,\,\frac{\sigma_{\varepsilon}^2}{1 - \varphi_2^2})$. This is a reasonable choice as an approximation because~\eqref{eq:ar2} and~\eqref{eq:ar1} share the same stationary distribution, $p(M_n)$, and the same one-step transition density, $p(M_{n+1}\mid M_n)$. Moreover, this choice aligns with a common approach to learning Markovian approximations from data, which is to try to match the statistics of single increments, i.e. to try to have the correct one-step transition density. However,~\eqref{eq:ar2} and~\eqref{eq:ar1} do, of course, differ. As mean-zero Gaussian processes they are uniquely determined by their autocovariance functions, $\gamma(m) = \mathrm{Cov}(M_n,\,M_{n+m})$. Let $\gamma$ denote the autocovariance of~\eqref{eq:ar2} and $\gamma_M$ denote that of~\eqref{eq:ar1}. Both satisfy linear recurrence relations (known as Yule--Walker equations~\citep{lutkepohl2013}),
\begin{subequations}\label{eq:gamma}
\begin{align}\label{eq:gammaa}
    \gamma(m) &= \varphi_1\gamma(m-1)+\varphi_2\gamma(m-2),\\
    \gamma_M(m) &= \varphi_M\gamma_M(m-1),\label{eq:gammab}
\end{align}
\end{subequations}
with solutions
\begin{subequations}\label{eq:gamma_explicit}
\begin{align}\label{eq:gammaa_explicit}
    \gamma(m) &= A \varphi_+^m + B \varphi_-^m,\\
    \gamma_M(m) &= \varphi_M^m \gamma(0),\label{eq:gammab_explicit}
\end{align}
\end{subequations}
where $\varphi_{\pm}=\frac{1}{2}\left(\varphi_1 \pm \sqrt{\varphi_1^2 + 4\varphi_2}\right)$ and $A$ and $B$ can be determined. Notably, for the parameter values we use, we have $\varphi_+=0.967$, $\varphi_-=-0.517$, and $\varphi_M=0.9$. Thus, for large $m$, $\gamma(m)$ is dominated by the contribution of the $\varphi_+^m$ term in~\eqref{eq:gammaa_explicit}, which in particular causes $\gamma(m)$ to decay slower than $\gamma_M(m)$. The AR$(2)$ process therefore exhibits longer-range memory than the natural AR$(1)$ approximation. On the other hand, with this in mind, we can use knowledge of $\varphi_+$ to propose a potentially better Markovian approximation: namely, the AR$(1)$ process
\begin{align}\label{eq:ar1+}
    M_n = \varphi_+ M_{n-1} + \epsilon ^{+}_n,
\end{align}
where $\epsilon ^{+}_n\overset{\text{iid}}{\sim}\mathcal{N}\left(0,\,\frac{\gamma(0)}{1 - \varphi_+^2}\right)$ and the variance of $\epsilon ^{+}_n$ is chosen to ensure this process has the same stationary distribution as~\eqref{eq:ar2} and~\eqref{eq:ar1}, namely
\begin{align}
    M_n\sim\mathcal{N}\left(0,\,\frac{(1 - \varphi_2)\sigma^2_{\varepsilon}}{(1+\varphi_2)(1-\varphi_1 - \varphi_2)(1+\varphi_1 - \varphi_2)}\right).
\end{align}

We thus consider non-Markovian model error modelled with two different Markovian parameterisations. With the first we aim to (crudely) represent standard Markovian parameterisations as they might be constructed for realistic models, albeit in a simplified setting where, importantly, the model error is independent of the model state. With the second, we highlight a potentially better Markovian parameterisation, which could only be identified from intimate knowledge of the model error process. We stress that in realistic modelling scenarios, where model error represents the effect of unresolved dynamics, the identification of better Markovian approximations is highly non-trivial.

We first consider differences in the stationary statistics of the Lorenz '63 system forced with these parameterisations. In particular we focus on the probability density function (pdf) of the first component of the system, and its temporal autocovariance function $r(\tau)$. For the pdf comparison we compute the Kullback--Leibler divergence and the Hellinger distance between the ``true'' system (that forced with AR$(2)$ noise) and (i) the unforced system, as in~\eqref{eq: l63}, (ii) the natural Markovian parameterisation, and (iii) the Markovian parameterisation with extended correlation, which we henceforth refer to as AR$(1)^{+}$. An intermediate density estimation step is performed using kernel density estimation to enable the calculation of these scores. To compare autocovariance functions we use the relative $L^2$ error
\begin{align}
    d_r(r,\,r') = \frac{\|r-r'\|_2}{\|r\|_2},
\end{align}
where each autocovariance function is computed up to time lags of $10$ Lyapunov times.
These scores, which we refer to as climate scores, are given in Table~\ref{tab:non_markovian_forcing}. Both parameterisations improve all three scores relative to the unforced system, but a much greater improvement is consistently seen with the AR$(1)^{+}$ parameterisation.

\begin{table}
    \centering
        \begin{tabular}{c|ccc}  
             & Unforced & AR$(1)$ & AR$(1)^+$  \\
             \hline
            KL & $2.1\times10^{-2}$ & $4.2\times10^{-3}$ & $9.8\times10^{-5}$ \\
            Hellinger & $7.1\times10^{-2}$ & $3.2\times10^{-2}$ & $4.9\times10^{-3}$ \\
            $d_r$ & $6.9\times10^{-1}$ & $2.4\times10^{-1}$ & $4.7\times10^{-2}$
        \end{tabular}
    \caption{Climate scores quantifying the discrepancy in stationary statistics of the Lorenz '63 system with non-Markovian AR$(2)$ forcing compared to Markovian parameterisations and the unforced system. Markovian parameterisations are seen to improve prediction, especially when chosen carefully to best approximate the non-Markovian error process.}
    \label{tab:non_markovian_forcing}
\end{table}

\begin{figure}
    \centering
    \includegraphics[width=5.5in]{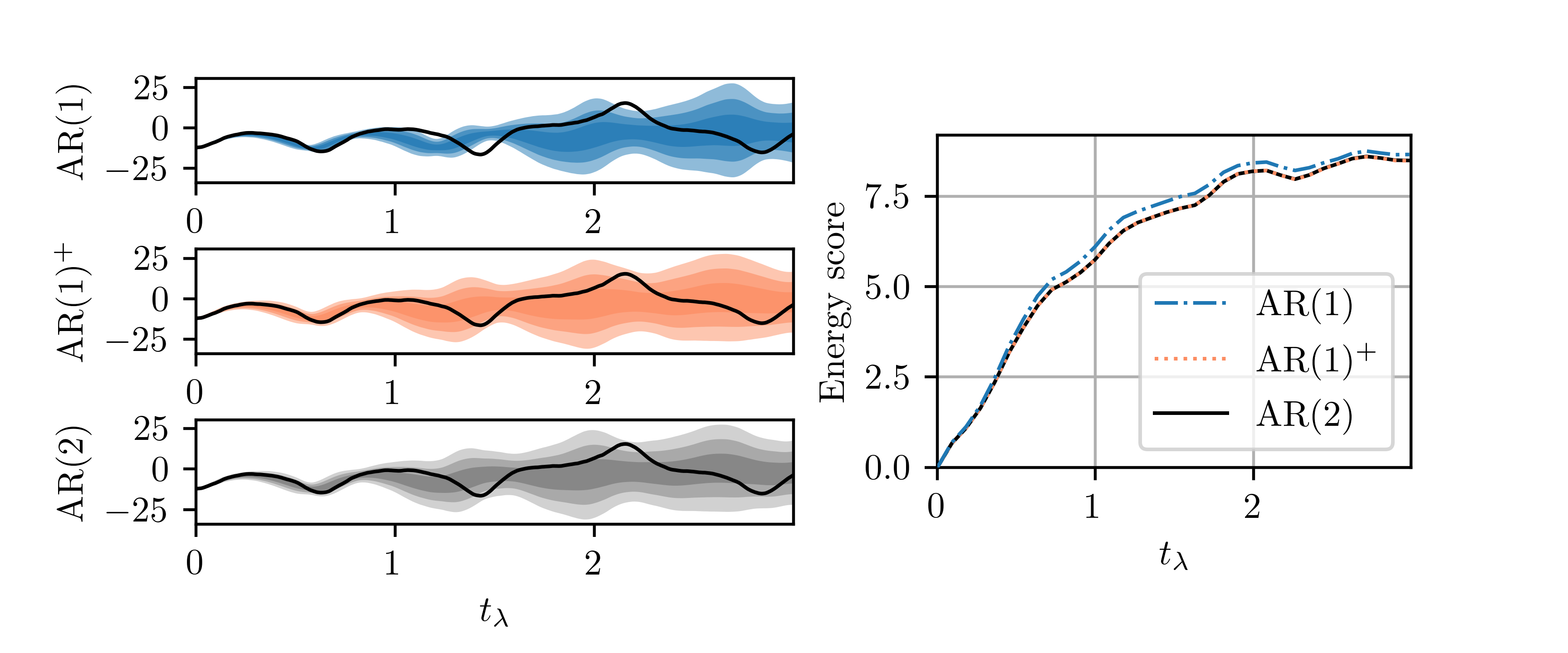}
    \caption{Left panel: ensemble forecasts of the Lorenz '63 system with Markovian $\mathrm{AR}(1)$ and $\mathrm{AR}(1)^{+}$ forcings and with non-Markovian $\mathrm{AR}(2)$ forcing. Black lines indicate a reference realisation of the system with $\mathrm{AR}(2)$ forcing and coloured regions indicate $\pm 3$ ensemble standard deviations about ensemble means. Right panel: mean energy score computed for the same parameterisations as a function of lead time.}
    \label{fig:l63_weather}
\end{figure}

We also assess the performance of these parameterisations in predicting at finite times, i.e. in the weather paradigm. To do so we generated $1000$ instances of ``weather'' from a single long trajectory of the true system, with each sufficiently separated in time to avoid noticeable correlation. We then simulated corresponding ensembles of size $100$ from each parameterised model from correct initial conditions. To assess the skill of these ensembles we use the energy score, a strictly proper scoring rule for ensemble forecasts~\citep{gneiting2007, gneiting2008}. For an ensemble of predictions $\{Z_1,\,\dots,\,Z_n\}$ of a variable $Z$ which takes a true value $Z^{\dagger}$, the energy score is
\begin{align}
    \mathcal{S}_{\mathrm{energy}}(\{Z_i\},\,Z^{\dagger}) = \frac{1}{n}\sum_{j=1}^n \|Z_j-Z^{\dagger}\| - \frac{1}{2n^2}\sum_{i=1}^n\sum_{j=1}^n\|Z_i - Z_j\|.
\end{align}
We compute the energy score, averaged over weather instances, as a function of lead time, for $t_{\lambda}\in[0, 3]$, where $t_{\lambda}$ is lead time normalised by the Lyapunov time of the unforced Lorenz '63 system. For comparison, we also generate ensembles with the true AR$(2)$ forcing, in order to determine optimal energy scores. The result is shown in Figure~\ref{fig:l63_weather}. We see that while the AR$(1)^{+}$ parameterisation leads to energy scores indistinguishable from the AR$(2)$ forcing, those attained with the AR$(1)$ parameterisation are noticeably worse across a range of lead times. For the sake of visualisation we also show in Figure~\ref{fig:l63_weather} a single weather instance and the ensembles produced with each forcing. As expected from the energy scores, the AR$(1)^{+}$ ensemble appears indistinguishable from the AR$(2)$ ensemble, but the AR$(1)$ ensemble shows noticeable error -- in particular, we see its predictions are overconfident, such that the true trajectory at times lies well outside the ensemble.

Overall, these results indicate that, in both the weather and climate paradigms, the performance of Markovian parameterisation in our setup is sensitive to the Markovian approximation used, and in particular, that the Markovian parameterisation that arises from correctly sampling one-step statistics is suboptimal. We suggest that this is likely the case also in stochastic parameterisation in Earth systems modelling broadly. It is therefore imperative that serious attention is given to the details of stochastic parameterisations implemented in large-scale modelling efforts, in particular with respect to the representation of memory.

\subsection{Spatially-correlated multiplicative VAR-1 forcing}\label{subsec: nonlocal}

Here we consider model error given by spatially-correlated multiplicative noise in the form $M_n = \Tilde{\bm{M}}_n \odot X_n$, driven by the VAR$(1)$ process
\begin{align}\label{eq:var1}
    \Tilde{\bm{M}}_n = \varphi \Tilde{\bm{M}}_{n-1} + \bm{\varepsilon}_n,\quad \bm{\varepsilon}_n\overset{\mathrm{iid}}{\sim}\mathcal{N}(\bm{0},\, \Sigma).
\end{align}
For clarity, note that $X_n$, $M_n$, $\tilde{\bm{M}}_n\in\mathbb{R}^3$ and $\odot$ denotes the Hadamard (or elementwise) product.
While it may seem more natural to continue with additive noise as in Section~\ref{subsec: non-markovian}, we found that spatially-correlated additive noise in the Lorenz '63 system led to a surprising phenomenon, where orbits lingered for extremely long times around one of the unstable steady states. We found that spatially-correlated multiplicative noise did not trigger this phenomenon, and so chose to avoid the issue, since we suspect it to be intrinsic to the Lorenz '63 system.

\begin{remark}
    We abuse the term space throughout. When we refer to spatial correlation or locality, we mean between components of $X^{(1)}_n$ or $M_n$, regardless of whether they correspond to the value of a single variable at various grid points in space, or entirely distinct physical variables. In particular, this means that the issues we discuss apply equally when the unresolved variables to be parameterised represent sub-grid scales of a variable being simulated (as in turbulence closure modelling) or variables which are neglected completely (e.g. biological components of the Earth system not featuring explicitly in a model).
\end{remark}

We set $\varphi=0.999$ such that the decorrelation time of $\Tilde{\bm{M}}_n$ is one Lyapunov time, and set
\begin{align}
    \Sigma = \kappa\begin{pmatrix}
        1 & \alpha & \alpha \\
        \alpha & 1 & \alpha \\
        \alpha & \alpha & 1
    \end{pmatrix},
\end{align}
with $\alpha=-0.45$ and $\kappa=1.81\times10^{-10}$ such that $\text{Var}(\tilde{\bm{M}}_n) = 10^{-7}$.

A spatially local parameterisation arises from neglecting spatial autocorrelation in $\bm{\varepsilon}_n$, i.e. by setting $\alpha=0$ so that $\Sigma=\kappa\mathbb{I}$. While in this simplistic setting there is negligible computational saving in adopting a spatially-local parameterisation, the incentive is clear in realistic applications, where modelling and sampling a high-dimensional conditional distribution $\bm{M}_n\mid\bm{X}_n$ is typically much more challenging and expensive than modelling each component as independent and identically (conditionally) distributed, i.e. considering $M_{n,\, i}\mid X_{n,\, i}$.

Table~\ref{tab:spatially_correlated_forcing} shows climate scores for the local parameterisation alongside those of the unforced system. Remarkably we see that the local parameterisation provides little improvement over the unforced system, indicating that assuming spatial locality can be significantly detrimental to stochastic parameterisations. On the other hand, we find that the local parameterisation performs well in the weather paradigm -- ensemble simulations and energy scores (not shown) were not significantly different from ones generated with the true spatially-correlated forcing.

The next section deals with the issues of non-Markovianity and spatial correlation in the case where model error arises, not synthetically from a chosen stochastic process, but from neglecting components of a dynamical system. The following experiments are, therefore, more representative of realistic modelling scenarios. In particular, the statistics of the error process are not known a priori, and stochastic parameterisations must be inferred from data.

\begin{table}
    \centering
        \begin{tabular}{c|cc}  
             & Unforced & Local VAR$(1)$ \\
             \hline
            KL & $2.7\times10^{-3}$ & $1.9\times10^{-3}$\\
            Hellinger & $2.6\times10^{-2}$ & $2.2\times10^{-2}$ \\
            $d_r$ & $4.4\times10^{-1}$ & $4.4\times10^{-1}$
        \end{tabular}
    \caption{Climate scores quantifying the discrepancy in stationary statistics of the Lorenz '63 system with spatially-correlated VAR$(1)$ forcing compared to a corresponding white-in-space parameterisation and the unforced system. The white-in-space parameterisation is seen to provide little improvement over the unforced system.}
    \label{tab:spatially_correlated_forcing}
\end{table}

\section{Experiments with the two-scale Lorenz '96 system}\label{sec: 96}

The experiments in this section make use of the Lorenz '96 system~\citep{lorenz1996}, usually written
\begin{subequations}
\label{eq: l96}
    \begin{align}
        \frac{\ud X_i}{\ud t} &= -X_{i-1}(X_{i-2}-X_{i+1}) - X_i + F - \frac{hc}{b}\sum_{j=J(i-1)+1}^{iJ}Y_j,\hspace{2em}i=1,\,\dots,\,I,\\
        \frac{\ud Y_j}{\ud t} &=-cbY_{j+1}(Y_{j+2} - Y_{j-1}) - cY_j - \frac{hc}{b}X_{\lfloor(j-1)/J\rfloor+1},\hspace{2em}j=1,\,\dots,\,IJ,
    \end{align}
\end{subequations}
with periodic boundary conditions, $X_{i+I}=X_i$ and $Y_{j+IJ}=Y_j$. The Lorenz '96 system has served for some time as a canonical testbed for research on parameterisation~\citep{wilks2005, crommelin2008, kwasniok2012, arnold2013, chorin2015, gagne2020, crommelin2021, levine2022, parthipan2023, bhouri2023} and data assimilation~\citep{hu2017, brajard2020, stanley2021} for several reasons: (i) it crudely mimics the character of chaotic advection--diffusion-type systems that appear in Earth systems modelling, (ii) it features two coupled variables $X$ and $Y$ and hence lends itself to the study of questions regarding model order reduction or partially observed dynamical systems, and (iii) it is inexpensive to simulate. The parameters $h$, $F$, $b$ and $c$ represent, respectively, a coupling strength coefficient, a forcing amplitude, a spatial scale ratio, and a time scale ratio. In other words $b$ and $c$ quantify scale separations between $X$ and $Y$, while $h$ quantifies the degree of influence each variable has on the other. In our experiments we fix parameters to canonical values $h=1$, $F=20$, $b=10$, $c=10$ with $I=8$ and $J=32$. In this regime $X$ can be thought of as a large scale variable and $Y$ an interacting small scale variable. A reduced order model for $X$ is obtained by neglecting entirely the term which couples $X$ to $Y$, giving
\begin{align}\label{eq: l96_reduced}
    \frac{\ud X_i}{\ud t} &= -X_{i-1}(X_{i-2}-X_{i+1}) - X_i + F,\hspace{2em}i=1,\,\dots,\,I.
\end{align}

To align notation with Section~\ref{sec: math} we once again write
\begin{align}
    X^{(1)}_{n+1} = \Psi_0\left(X_n^{(1)}\right) + M_n,
\end{align}
where $X^{(1)}=X$, $\Psi_0$ denotes the map of the reduced system~\eqref{eq: l96_reduced} (after having been discretised in time), and $M_n$ represents model error relative to the full dynamics of~\eqref{eq: l96}. Thus, in this section the reference ``true'' dynamics is that of~\eqref{eq: l96} with an imperfect model given by~\eqref{eq: l96_reduced}.

\subsection{Data-driven parameterisations}
\subsubsection{Data}
We consider parameterisations learned from data. In particular, we 
simulate~\eqref{eq: l96} using fourth-order Runge--Kutta timestepping with stepsize $\Delta t=10^{-3}$. Then, by evaluating $\Psi_0$ at each value of $X^{(1)}$ observed in simulation, we diagnose corresponding values of $M_n$. A climate dataset is generated by simulating for $10^4$ model time units from a random initial condition ($X_i(0),\,Y_j(0) \overset{\text{iid}}{\sim}\mathcal{N}(0,\,1)$), giving
\begin{align}
    \mathcal{D}_{\text{climate}}\coloneqq \left\{(X^{(1)}_n,\,M_n)_{n=0}^{10^7}\right\}.
\end{align}
A corresponding weather dataset is constructed by partitioning the climate dataset into $N_w=10^3$ slices, each of $L=10^4$ timesteps, corresponding to $10$ model time units, giving
\begin{align}
    \mathcal{D}_{\text{weather}}\coloneqq \left\{\left\{(X^{(1)}_n,\,M_n)_{n=L(k-1)}^{Lk}\right\},\, i = 1,\,\dots,\, N_w \right\}.
\end{align}
We use the climate dataset to fit parameterisations and later to evaluate how well each parameterisation reproduces climate statistics. The weather dataset is used to evaluate how well each parameterisation represents uncertainty in weather prediction.

\subsubsection{Parameterisations}
We construct Markovian parameterisations, using probabilistic neural networks (specifically mixture density networks (MDNs)~\citep{bishop1994}) to model the conditional distribution $M_n\mid X^{(1)}_n$. MDNs combine a neural network with a Gaussian mixture model to produce a highly flexible parametric form for conditional densities. In this case, this means, for each possible value of $X^{(1)}_n$, the conditional density $p(M_n\mid X^{(1)}_n)$ is modelled by the density of a Gaussian mixture distribution, whose parameters are given as a function of $X^{(1)}_n$ modelled with a neural network. The neural network is fit using maximum likelihood estimation. A detailed description can be found in~\citet{brolly2023} where MDNs are used to model the transition density of ocean surface drifters. Here we consider Gaussian mixtures with $32$ components and neural networks with $4$ hidden layers, each having $128$ neurons and $\tanh$ activation function. The parametric form of the Gaussian mixture is chosen to be flexible and, in particular, allow non-Gaussianity to be captured. Note that the number of parameters required to specify a Gaussian (mixture) distribution for $M_n$, and hence, the dimension of the neural network's output, is $\mathcal{O}(d_1^2)$. Since $d_1=8$ in this experiment, this is not a problem, but in realistic applications where $d_1$ is large, other methods may be required to more efficiently model $p(M_n\mid X^{(1)}_n)$. The procedure of simulating with Markovian MDN parameterisations is outlined in Algorithm~\ref{algo: mdn_sim}, wherein $\bm{f}_{\text{NN}}$ denotes the trained neural network and $\bm{\theta}$ denotes the parameters of the Gaussian mixture.

\begin{algorithm}
    \caption{Simulating with Markovian MDN parameterisations.}\label{algo: mdn_sim}
    \begin{algorithmic}[1]
    \State $X^{(1)}_0 \gets x_0$
    \For{$n = 0$ to $N-1$}
        \State $X^{(1)}_{n+1} \gets \Psi_0(X^{(1)}_n)$
        \State $\bm{\theta} \gets \bm{f}_{\text{NN}}(X^{(1)}_n)$\\
        \hspace{1.5em}draw $M_n \sim \text{GaussianMixture}(\bm{\theta})$
        \State $X^{(1)}_{n+1} \gets X^{(1)}_{n+1} + M_n$
    \EndFor\\
    \Return $X^{(1)}_0,\,\cdots,\,X^{(1)}_{N}$
    \end{algorithmic}
\end{algorithm}

\begin{remark}
    The choice to model $M_n\mid X^{(1)}_n$ rather than $M_n\mid X^{(1)}_n,\,M_{n-1}$ is a purely pragmatic one. As others have noted~\citep{gagne2020, parthipan2023}, there is a tendency for data-driven models of $M_n\mid X^{(1)}_n,\,M_{n-1}$ to ignore information from $X^{(1)}_n$, because of high correlation between $M_n$ and $M_{n-1}$, leading to poor performance.
\end{remark}


Spatially local parameterisations can be constructed similarly while enforcing constraints on the form of $p(M_n\mid X^{(1)}_n)$. We distinguish two forms of spatial locality. We call a parameterisation weakly local if the components of $M_n$ are modelled as conditionally independent given $X^{(1)}_n$, i.e. $p(M_n\mid X^{(1)}_n) = \prod_i p(M_{n,\, i}\mid X^{(1)}_n)$. In the MDN approach this can be enforced by assuming diagonal covariance matrices for the Gaussian mixture components. Note that the number of parameters required to specify a Gaussian (mixture) distribution is then $\mathcal{O}(d_1)$. We call a parameterisation strongly local if additionally the components of the error, $M_{n,\, i_1}$, are independent of $X^{(1)}_{n,\,i_2}$ for $i_2\neq i_1$, i.e. $p(M_n\mid X^{(1)}_n) = \prod_i p(M_{n,\, i}\mid X^{(1)}_{n,\,i})$. In the strongly local case the number of parameters required to specify a Gaussian (mixture) distribution is independent of $d_1$. Previous studies on stochastic parameterisation in the Lorenz '96 system have predominantly adopted the strongly local assumption~\citep{crommelin2008, kwasniok2012, arnold2013, chorin2015, gagne2020, parthipan2023}. The appeal of spatial locality clearly lies in reduced computational complexity, both in the construction and in the deployment of parameterisations. But does spatial locality affect the performance of parameterisations? Here we deploy both local and nonlocal parameterisations and compare them on the basis of climate and weather scores.

We restrict attention in this section to Markovian parameterisations. However, considering the experiments of Section~\ref{sec: 63}, we should expect that assuming Markovianity may be detrimental. \citet{parthipan2023} implemented a non-Markovian (though strongly local) parameterisation in the Lorenz '96 system. Their parameterisation was based on a recurrent neural network and they saw improvement in performance over a baseline Markovian model. However, it remains to be seen whether non-Markovian parameterisations of this type can be successfully implemented in high-dimensional models. In Section~\ref{sec: 63} we noted that the Markovian parameterisation that arises from correctly modelling the one-step transition density is not in general the optimal Markovian parameterisation. Thus, as a compromise we ask: how can we improve on current Markovian models?

As a simple modification of the Markovian parameterisations described above, we propose to introduce a parameter $t_p\in\mathbb{Z}_{>0}$. And rather than sample $M_n$ at each timestep, we repeatedly hold $M_n$ fixed for $t_p$ timesteps at a time before sampling again. We call $t_p$ the timestep of the parameterisation. Intuitively, setting $t_p>1$ induces increased memory in the error process, in a similar manner to the AR$(1)^+$ parameterisation used in Section~\ref{sec: 63}. An attractive feature of this modification is that it actually reduces the cost of parameterisations. Given that the dominant cost in deploying these parameterisations is sampling the conditional distribution $M_n\mid X^{(1)}_n$, we can expect a reduction in cost of a factor of $t_p$. More importantly, by optimising $t_p$ to maximise performance, either with respect to climate or weather scores, as appropriate for the application, we may be able to improve Markovian parameterisations easily. We note that a similar approach was followed by~\citet{zanna2017scale} in constructing stochastic parameterisations of sub-grid-scale forcing in a quasi-geostrophic ocean model.

\subsubsection{Results}

\begin{figure}
    \centering
    \includegraphics[height=8in]{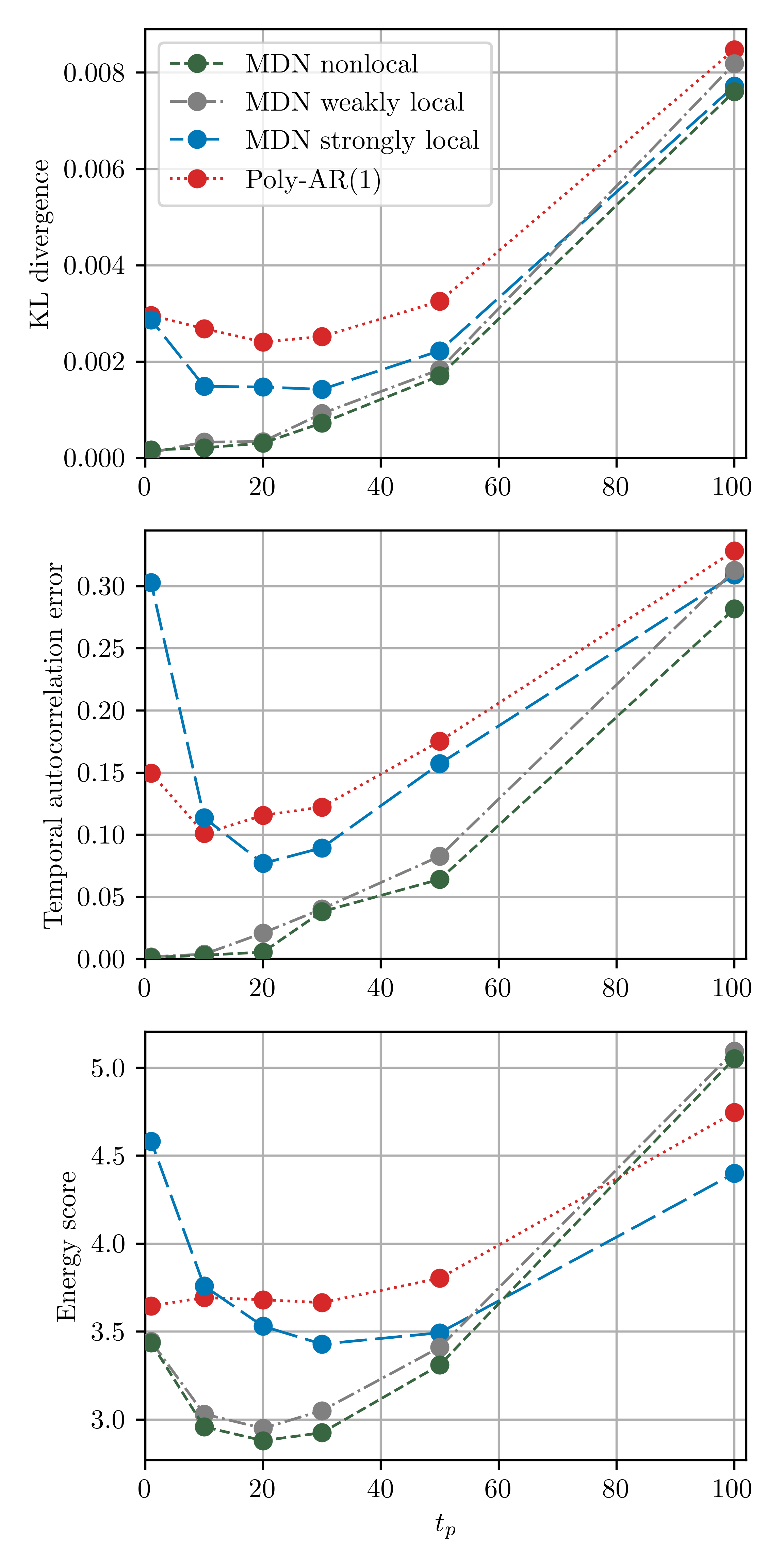}
    \caption{KL divergence between the true stationary measure of $X_k$ in the Lorenz '96 model (upper panel), error in the temporal autocorrelation of $X_k$ (middle panel), and energy score at time $t=1$ (lower panel), in the parameterised Lorenz '96 models as a function of parameterisation time step.}
    \label{fig:l96_scores}
\end{figure}

Figure~\ref{fig:l96_scores} shows climate and weather scores for nonlocal, weakly local, and strongly local Markovian MDN parameterisations for $t_p\in\{1,\,10,\,20,\,30,\,50,\,100\}$. For the temporal autocorrelation error, autocorrelations functions are computed up to time lags of $10$ model time units. Scores are also shown for a baseline parameterisation, named Poly-AR$(1)$, used by~\citet{arnold2013} in the Lorenz '96 system. Poly-AR$(1)$ models $M_n$ using a polynomial function of $X^{(1)}_n$ for the conditional mean $\mathbb{E}[M_n\mid X^{(1)}_n]$ and an AR-$(1)$ process for the residuals $M_n - \mathbb{E}[M_n\mid X^{(1)}_n]$; it is Markovian, strongly local, and Gaussian. For each parameterisation, scores vary strongly with $t_p$. Comparing parameterisations by the scores they achieve with their optimal value of $t_p$, there is a clear ordering of the parameterisations which persists across all of the scores considered: the nonlocal MDN parameterisation outperforms the weakly local MDN parameterisation, which outperforms the strongly local MDN parameterisation, which outperforms the strongly local and Gaussian Poly-AR$(1)$ parameterisation. This ordering reflects an intuitive conclusion that parameterisations which rely on fewer simplifying assumptions are quantifiably better.

Tuning $t_p$ is beneficial if and only if the optimal value is greater than one. In most cases we see that this is the case, but in the case of the nonlocal and weakly local MDN parameterisations we see that climate scores are optimised by taking $t_p=1$. On the other hand the optimal value for the weather score (the energy score) for those same parameterisations is around $t_p=20$, highlighting that parameterisations can (and perhaps should) be tailored towards the particular needs of the modeller. Moreover, the observed difference in optimal values between weather and climate scores suggests that the performance of parameterisations in weather forecasting is not necessarily indicative of performance in climate modelling.

In terms of the climate scores, we see the best performance with the MDN nonlocal and weakly local parameterisations with $t_p=1$. Noting that these two parameterisations perform indistinguishably with respect to these scores, we also constructed an additional nonlocal but deterministic parameterisation, having the same architecture as the MDN but outputting only a point estimate of $\bm{M}_n$ rather than $\bm{\theta}$, and trained to minimise a standard mean squared error loss. This deterministic parameterisation was found to score equally well in terms of KL divergence and temporal autocorrelation error. This suggests that, in some cases, deterministic parameterisations may be sufficient for climate modelling purposes, provided that they are nonlocal. However, this may not generalise to other models, given that the configuration of the Lorenz '96 model considered here features significant scale separation that may be absent in many applications.

\begin{figure}
    \centering
    \includegraphics[width=5.5in]{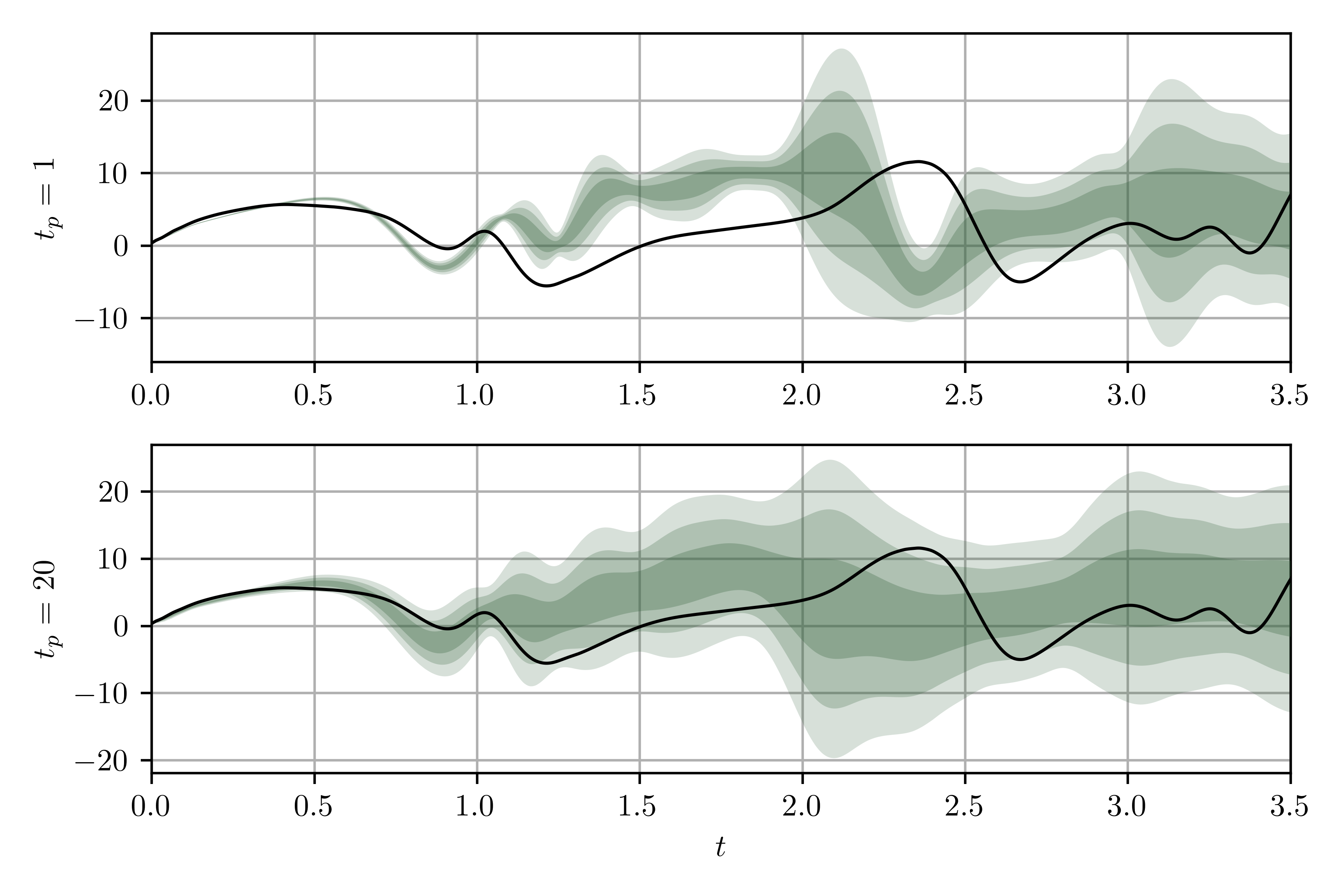}
    \caption{Ensemble forecasts from the nonlocal MDN-parameterised Lorenz '96 model with $t_p=1$ (upper panel) and $t_p=20$ (lower panel). The black line indicates the true trajectory of $X_1$. Contours show $\pm3$ ensemble standard deviations about the ensemble mean.}
    \label{fig:l96_ensemble_tp_comp}
\end{figure}

In Figure~\ref{fig:l96_ensemble_tp_comp} we show how ensembles simulated with the nonlocal MDN parameterisation are affected by tuning $t_p$. Ensembles are shown for $t_p=1$ and for the value that optimises the energy score, $t_p=20$. With $t_p=1$ the ensemble is significantly biased and overconfident, so that the true trajectory generated with the full Lorenz '96 system lies well outside the spread of the ensemble. With $t_p=20$ this behaviour is not seen, and the ensemble appears to represent uncertainty much more plausibly. While this is quantified better by the energy scores reported above, Figure~\ref{fig:l96_ensemble_tp_comp} provides a visual indication of the improvement that tuning $t_p$ can bring in the weather paradigm. Plots of this kind are notably absent from the literature on stochastic parameterisation. 
We recommend their use as a simple means of identifying deficiencies in uncertainty quantification provided by stochastic parameterisations. If reference trajectories lies well outside the range of ensembles, there is clear cause for concern. The only complication in producing  plots of this type in operational forecasting situations is the choice of a representative scalar observable; in Figure~\ref{fig:l96_ensemble_tp_comp} the first component of the Lorenz '96 model is shown, somewhat arbitrarily. In realistic Earth system models one can plot, for example, the value of an important prognostic variable at a certain point in space, or the average over a region. However, since performance viewed through these plots may be sensitive to the choice of observable, modellers should take care and consider them alongside proper scoring rules for multivariate ensembles.

While the technique of tuning $t_p$ is far from a final solution to the problem of finding optimal Markovian parameterisations, it could be used readily in situations where making major changes to existing parameterisations is difficult. In operational settings one could perform a small number of experimental short-term forecasts where $t_p$ is varied over a range of values and inspect plots similar to Figure~\ref{fig:l96_ensemble_tp_comp}. The energy score could be used to make the results quantitative.

\section{Conclusion}\label{sec: conc}

In this work we highlight the impact of locality assumptions on the performance of stochastic parameterisations. We show that, even in simple settings, assuming Markovianity and/or spatial locality can be detrimental to model predictions, in both the weather and climate paradigms. On the other hand we show that the impact of assuming Markovianity can be lessened by careful choice of Markovian parameterisation. In particular, we show that Markovian parameterisations which correctly represent the statistics of one-step transitions are generally suboptimal. We introduce a simple modification which can be made to Markovian parameterisations of this type, whereby tuning a single parameter can yield significant improvements in predictive performance. Finally, in tuning our parameterisations, we observe a divergence between configurations which perform best in predicting weather and those which best represent climate --- a reminder that, so long as parameterisations remain imperfect, optimal performance in short-term prediction is no guarantee of optimal performance in modelling long-term behaviour, and vice versa.
\vspace{1em}

\noindent {\large\textbf{Acknowledgements}}

\normalsize
\noindent I thank Aretha Teckentrup for her support during the preparation of this article. I also thank James R. Maddison and Jacques Vanneste for helpful feedback on an early version of the manuscript. Finally, I thank the three anonymous reviewers for their constructive comments. This work was supported by EPSRC grant number EP/X01259X/1.

\bibliography{main}

\end{document}